\documentclass{article}
\usepackage{latexsym}
\usepackage[all]{xy}
\usepackage{amssymb,amsmath}
\def\cur{\Gamma}
\def\Pic{\mathrm{Pic}}
\def\Hom{\mathrm{Hom}}
\def\Supp{\mathrm{Supp}}

\def\P{\mathbb{P}}
\def\cP{\mathcal{P}}

\def\O{\mathcal{O}}

\def\X{\mathcal{X}}

\def\I{\mathcal{I}}

\def\M{\mathcal{M}}

\def\dim{\mathrm{dim}}

\def\separation{\medskip}
  \newtheorem{theorem}{Theorem}[section]
\newtheorem{lemma}[theorem]{Lemma}
\newtheorem{prop}[theorem]{Proposition}
\newtheorem{definition}[theorem]{Definition}
\newtheorem{corollary}[theorem]{Corollary}
\newtheorem{remarks}[theorem]{Remarks}
\newtheorem{remark}[theorem]{Remark}

\newcommand{\proof}{{\it Proof.}\ }
\newcommand{\qed}{\hfill  $\Box$\separation}

\begin{document}

\title{The symmetric square of a curve \\ and  the Petri map}
 \author{A. BRUNO - E. SERNESI \footnote{Both authors are members of GNSAGA-INDAM}}
\date{}
\maketitle
\abstract{Let $\M_g$ be the course moduli space
of complex projective nonsingular curves of genus $g$. We prove
that when the Brill-Noether number $\rho(g,1,n)$ is non-negative
the Petri locus $P^1_{g,n}\subset \M_g$ has a divisorial component whose closure has a non-empty intersection with $\Delta_0$. In order to prove
the result we show that the scheme $G^1_n(\Gamma)$ that parametrizes degree $n$ pencils on a curve $\Gamma$ is isomorphic to a component of the Hilbert scheme parametrizing certain curves on the symmetric square $\Gamma_2$ of $\Gamma$ and we study the properties of such a family of curves.}

\section{Introduction}

This paper is devoted to  some aspects of the Petri condition for
linear pencils on a complex algebraic curve.  Let $C$ be a
nonsingular irreducible projective curve of genus $g \geq 2$
defined over $\mathbb{C}$, and let $(L,V)$ be a pair
 consisting of an invertible sheaf  $L$ on $C$ and of an $(r+1)$-dimensional
 vector  subspace $V \subset H^0(L)$,  $r \ge 0$. $(L,V)$   is called
  a \emph{linear series of dimension $r$ and degree $n$}, or a $g^r_n$.
  If $V=H^0(L)$
 then the $g^r_n$ is said to be \emph{complete}.
A \emph{linear pencil} is a $g^1_n$.

\noindent If $(L,V)$ is a $g^r_n$ then the \emph{Petri map}  for
$(L,V)$ is the natural multiplication map
\[
\xymatrix{ \mu_0(L,V):V \otimes H^0(\omega_C L^{-1}) \ar[r]&
H^0(\omega_C)}
\]
The Petri map for $L$ is
\[
\xymatrix{ \mu_0(L):H^0(L) \otimes H^0(\omega_C L^{-1}) \ar[r]&
H^0(\omega_C)}
\]
Recall that $C$ is called a \emph{Petri curve} if the Petri map
$\mu_0(L)$ is injective for every invertible sheaf $L$ on $C$. We
will say that  $C$ is \emph{Petri with respect to $g^r_n$'s} if
the Petri map $ \mu_0(L,V)$
   is injective for every $g^r_n$ $(L,V)$ on $C$.
    We say that \emph{$C$ is Petri with respect to pencils} if it is Petri w.r. to $g^1_n$'s for every $n$.

By the Gieseker-Petri theorem \cite{dG82}
 we know that in $\M_g$, the course moduli space of nonsingular projective curves of genus $g$,
  the locus of curves which are not Petri  is a proper closed subset $P_g$,
  called the \emph{Petri locus}.
This locus decomposes as
\[
P_g = \bigcup_{r,n} P^r_{g,n}
\]
where we  denote by  $P^r_{g,n}\subset \M_g$ the locus of curves
which are not Petri w.r. to $g^r_n$'s. The structure of $P_g$ and
of  $P^r_{g,n}$ is not known in general: they might a priori have
several components and not be equidimensional. If the
Brill-Noether number
\[
\rho(g,r,n) := g -(r+1)(g-n+r)
\]
is nonnegative then  it is conjectured that $P^r_{g,n}$ has pure
codimension one if it is non-empty. This is known to be true  in
some special cases, namely when $\rho(g,r,n)=0$,
 and  for $r=1$ and $n=g-1$ \cite{mT88}.
 In \cite{BS11} we proved that if
  $\rho(g,r+1,n)<0$   then $P^r_{g,n}$ has pure
codimension one outside $P^{r+1}_{g,n}$.
 When $g \le 13$ then $P_g$ has pure codimension one \cite{aC05,aC07,LC11}.

 Denote by $\overline{\M}_g$ the moduli space of stable curves, and let
\[
\overline{\M}_g \backslash \M_g =\Delta_0\cup \cdots \cup
\Delta_{[{g\over 2}]}
 \]
 be its boundary, in standard notation.
 In \cite{gF05} G.  Farkas  has   proved the existence of at least one
 divisorial component of $P^1_{g,n}$     in case $\rho(g,1,n) \ge 0$ and $n \le g-1$,
 using the theory of limit linear series.  He found a divisorial component which has a
 nonempty intersection with  $\Delta_1$. Another proof has been given in \cite{CT08},
  by degeneration to a stable curve with $g$ elliptic tails. The
  method of \cite{gF05} has been extended in \cite{gF08} to
  arbitrary $r$.

  Our main result is:

\begin{theorem}\label{T:main2}
Assume that $g,n \ge 3$    are such that  $\rho(g,1,n) \ge 0$.
Then $P^1_{g,n}$ has a divisorial component whose closure has a
non-empty intersection with $\Delta_0$.
\end{theorem}

We obtain this result by induction on $g$, via a preliminary study
of independent interest  of certain pencils on curves of genus
$g-1$. Namely, we consider a nonsingular curve $\cur$ of genus
$\gamma= g-1$ which is Petri w.r. to pencils, and, for a given
$x+y \in \cur_2$, the second symmetric power of $\cur$, we
consider those $g^1_n$'s $(L,V)$ such that $x+y$ is contained in
one of their  divisors: we call them \emph{neutral linear pencils}
w.r. to $x+y$. They form a divisor $N^1_n(x+y)$ of $G^1_n(\cur)$
which is closely related to the set of \emph{generalized pencils}
 on the nodal curve $X$ of arithmetic genus $g$ obtained from
$\cur$ after identifying $x=y$ (see Definition \ref{D:cusp1}). We
undertake a detailed study of the divisors $N^1_n(x+y)$, of their
singular locus, and of their relation with generalized pencils on
$X$. We  prove that the only generalized $g^1_n$'s on $X$ which
are non-Petri are those corresponding to singular elements of
$N^1_n(x+y)$ for which neither $x$ nor $y$ are base points. Then
we show that those $x+y \in \cur_2$ for which $N^1_n(x+y)$
contains such singularities  are contained in a curve $R_n \subset
\cur_2$, and that in such a case the singularities of that sort
are finitely many (Proposition \ref{P:NP2}). In particular, by
choosing $x+y$ outside all such curves for the possible (finite)
set of values of $n$, the curve $X$ will be Petri w.r. to
generalized pencils. From this result \emph{we obtain an inductive
proof of the Gieseker-Petri theorem for pencils.} The proof of
Proposition \ref{P:NP2} uses Proposition \ref{P:kermu} asserting
that in $\cur$ every invertible sheaf with $h^0(L) \ge 2$
satisfies $H^1(L^2(-x))=0$ for all $x \in \cur$.  This statement
has far-reaching consequences and has been the starting point of
this work.

A key role here is
  played by   the curves
$s_{L,V}\subset\cur_2$ associated to  $g^1_n$'s $(L,V)$,
    consisting of all $x+y \in \cur_2$ which are contained in some
    divisor of $(L,V)$. The family of all such curves is studied
    in detail: in particular we show that the scheme $G^1_n(\cur)$
    of $g^1_n$'s, which parametrizes the curves $s_{L,V}$,
    is isomorphic to a component of the Hilbert scheme of
     $\cur_2$  (Theorem \ref{T:hilb}).
A closer look at the above mentioned curve   $R_n \subset \cur_2$
shows that it does not contain the diagonal. For this part we use
the main result of \cite{EH83}, which implies that   $\Delta_1
\subset \overline{\M}_g$ is not contained in the closure of $P_g$.
We also need to prove that $R_n$ is non-empty. Here we apply a
theorem of Steffen \cite{fS98} to the characteristic map of the
the family of $s_{L,V}$'s (Proposition \ref{P:NP1}).

With the aid of all the previous results we
   prove   Theorem \ref{T:main2} as follows. We fix a non-Petri nodal curve $X$ as above
 and we consider a modular family
   of stable genus $g$ curves $f:\X \longrightarrow B$ containing
   $X$  among its fibres. A straightforward construction
   produces   an irreducible
   variety $\widetilde{P}$ of dimension $\ge 3g-4$ parametrizing pairs $(C,(L,V))$ such that $(L,V)$ is non-Petri and a
    morphism $\widetilde{P}\longrightarrow B$.  Since this map is
    finite above $X$ it follows that its image is a divisorial
    component $P$ of $P^1_{g,n}$. This component contains $X$ and
    therefore it   intersects $\Delta_0$.

    The structure of the paper is as follows.
In \S \ref{S:bpfpt}  we discuss the properties of the Petri map
for pencils.  In \S \ref{S:symm} we recall some cohomological
properties of line bundles on $\cur_2$ and we introduce and study
the curves $s_{L,V}$. The schemes of neutral linear series are
introduced in general in \S \ref{S:NLS} and they are studied in
detail in \S  \ref{S:NLPII} in the special case of pencils. In \S
\ref{S:singp1} we introduce and study the curves $R_n$. Then in \S
\ref{S:nodcur} we interpret in terms of nodal curves  what has
been proved before. The final \S \ref{S:main2} is devoted to the
proof of Theorem \ref{T:main2}.

\section{The Petri map for pencils}\label{S:bpfpt}

Let $\cur$ be a projective nonsingular irreducible curve of genus
$\gamma \ge 3$,
 $L$ a line bundle on $\cur$ and $V \subset H^0(L)$ a subspace.
 If $D$ is an effective divisor on $\cur$ we define
\[
V(-D) :=  V \cap H^0(L(-D))
\]
and \[ V(D) := j(V) \]
 where $j:H^0(L) \longrightarrow H^0(L(D))$
is the inclusion.

 \emph{Assume that $(L,V)$ is a  base-point free pencil} on $\cur$, i.e. a linear
 pencil such that $V$  generates $L$. The Koszul complex of $V$ gives  the exact sequence:
\begin{equation}\label{E:LV}
\xymatrix{ 0 \ar[r] & \bigwedge^2V\otimes L^{-2} \ar[r] & V
\otimes L^{-1} \ar[r] & \O_\cur \ar[r] & 0 }
\end{equation}
If we tensor it by $LM$ where $M$ is any line bundle, we obtain:
\[
\xymatrix{ 0 \ar[r] &  \bigwedge^2V\otimes L^{-1}M \ar[r] & V
\otimes M \ar[r] & LM \ar[r] & 0 }\] and taking global sections we
obtain the exact sequence:
\begin{equation}\label{E:LM1}
\xymatrix{
 0\ar[r]& \bigwedge^2V\otimes H^0(L^{-1}M)\ar[r]&V \otimes H^0(M) \ar[r]^-m & H^0(LM)}
\end{equation}
where  $m$ is the multiplication map. Therefore we have a
canonical identification:
\begin{equation}\label{E:LM2}
\ker(m)= \bigwedge^2V\otimes H^0(L^{-1}M)
\end{equation}
If we take $M=\omega_\cur L^{-1}$  then   (\ref{E:LM1}) becomes
\begin{equation}\label{E:LV1}
\xymatrix{
 0\ar[r]& \bigwedge^2V\otimes H^0(\omega_\cur L^{-2})\ar[r]&
 V \otimes H^0(\omega_\cur L^{-1}) \ar[r]^-{\mu_0} & H^0(\omega_\cur)}
\end{equation}
where  $\mu_0=\mu_0(L,V)$. Therefore we have a canonical
identification:
\[
\ker(\mu_0(L,V))= \bigwedge^2V\otimes H^0(\omega_\cur L^{-2})
\]
After tensoring (\ref{E:LV}) by $L^2$   we have an exact sequence:
\begin{equation}\label{E:LV2}
\xymatrix{ 0 \ar[r] &  \bigwedge^2V\otimes\O_\cur \ar[r] & V
\otimes L \ar[r] & L^2 \ar[r] & 0 }
\end{equation}
which induces:
\begin{equation}\label{E:rho}
\xymatrix{
 \bigwedge^2V\otimes H^1(\O_\cur) \ar[r]^-{\rho^\vee} & V \otimes H^1(L) \ar[r] & H^1(L^2)\ar[r] & 0 }
\end{equation}
and dualizing:
\[
\xymatrix{ 0\ar[r]&H^0(\omega_\cur L^{-2}) \ar[r] &V^\vee\otimes
H^0(\omega_\cur L^{-1}) \ar[r]^-\rho &
 \bigwedge^2V^\vee \otimes H^0(\omega_\cur)}
\]
If we now tensor by $\bigwedge^2V$ and recall that
$V=\bigwedge^2V\otimes V^\vee$ we obtain the sequence
(\ref{E:LV1}) again. This implies that $\mu_0(L,V)$ and $\rho$ are
the same map. Therefore we have  an isomorphism
\begin{equation}\label{E:kermu}
\ker(\mu_0(L,V)) \cong \ker(\rho) = H^1(L^2)^\vee
\end{equation}

   \emph{Suppose now that $(L,V)$ is any linear pencil on $\cur$} and
   let $B$ be its fixed divisor. Then, applying the above analysis to
   the fixed point free linear pencil $(L(-B),V(-B))$ we obtain that, for any line bundle $M$,
  \begin{equation}\label{E:LM3}
  \begin{array}{ll}
  \xymatrix{
  \mathrm{ker}[V\otimes H^0(M)\ar[r]^-m&H^0(LM)] }\cong \\
  \xymatrix{
  \mathrm{ker}[V(-B)\otimes H^0(M)\ar[r]&H^0(LM(-B))]}
  &\cong H^0(L^{-1}M(B))
  \end{array}
  \end{equation}
  In particular we obtain:
  \[
 \ker[\mu_0(L,V)] \cong \xymatrix{\ker[V(-B)\otimes H^0(\omega_\cur L^{-1})\ar[r]& H^0(\omega_\cur(-B))]} \cong
    H^0(\omega_\cur L^{-2}(B))
 \]
 We have the following proposition:

  \begin{prop}\label{P:kermu}
  If $\cur$ is Petri with respect to pencils then:
  \begin{itemize}

  \item[(i)]  $H^1(L^2)=0$ for every $L\in\Pic(\cur)$ such that $h^0(L)\ge 2$.

  \item[(ii)] $H^1(L^2(-x))=0$ for every $L\in\Pic(\cur)$ such that $h^0(L)\ge 2$ and for every $x \in \cur$.
  \end{itemize}
  \end{prop}

  \proof We only need to prove (ii). Let $L \in \Pic(\cur)$ be such that $h^0(L) \ge 2$ and
  assume that $H^1(L^2(-x)) \ne 0$ for some $x \in \cur$. Consider the exact sequence:
  \[
  \xymatrix{
  0\ar[r]& L^2(-x) \ar[r]&L^2 \ar[r] & L^2\otimes\O_x \ar[r]& 0 }
  \]
If $x$ is not base-point of $L^2$ then $h^1(L^2) > 0$ and this
contradicts (i).
 If $x$ is base-point of $L^2$ then it is also a base point of $L$. Therefore
  $L=M(x)$ with $h^0(M)=h^0(L)\ge 2$. We have $L^2(-x)=M^2(x)$ so that:
\[
0 < h^1(L^2(-x)) = h^1(M^2(x)) \le h^1(M^2)
\]
again contradicting (i). \qed

A special case of part (ii) of the proposition is in \cite{gP03},
Lemma 4.2.

Some of these remarks can be extended to linear series of higher
dimension, as follows.
 If $\dim(V)=r+1 \ge 3$ and $V$ generates $L$, from   the Koszul complex we obtain the exact sequence:
\[
\xymatrix{ 0\ar[r]&\bigwedge^{r+1}V\otimes L^{-r-1}\ar[r] &
\bigwedge^rV\otimes L^{-r}\ar[r]&
\varphi_V^*\Omega^{r-1}_{\P}\ar[r]&0 }
\]
where $\varphi_V: \cur \to \P(V^\vee)=\P$ is the morphism defined
by $V$. After tensoring by $L^{r+1}$ and taking cohomology we
obtain the exact sequence:
\[
\xymatrix{ \bigwedge^{r+1}V\otimes H^1(\O_\cur)
\ar[r]^-{\rho^\vee} & \bigwedge^rV\otimes H^1(L) \ar[r]&
H^1(\varphi_V^*\Omega^{r-1}_{\P}(r+1)) \ar[r] & 0}
\]
Again, the map $\rho^\vee$ can be identified with
$\mu_0(L,V)^\vee$, after having identified
\[
V^\vee = \bigwedge^{r+1}V^\vee\otimes \bigwedge^rV
\]
and   we find the well-known isomorphism:
\[
\ker(\mu_0(L,V)) \cong
H^1(\cur,\varphi_V^*\Omega^{r-1}_{\P}(r+1))^\vee
\]

\section{The second symmetric power}\label{S:symm}

 Let $\cur_2$ be the second symmetric power of $\cur$.
Denote by
\[
\xymatrix{ \sigma: \cur\times\cur \ar[r] & \cur_2}
\]
the natural degree two quotient morphism and by
\[
\xymatrix{ \cur&\cur\times\cur \ar[l]_{\pi_1} \ar[r]^{\pi_2} &
\cur }
\]
the projections. Given $L,M$ invertible sheaves on $\cur$ we will
let
\[
L\boxtimes M := \pi_1^*L \otimes \pi_2^* M
\]
 Let $D \subset \cur\times\cur$ be the diagonal and
\[
\delta = \sigma(D) = \{2x: x \in \cur\}\subset \cur_2
\]
We have
\[
\sigma^*(\delta) = 2D
\]
Moreover, the invertible sheaf $\O_{\cur_2}(\delta)$, which we
will shortly denote by $\delta$, is divisible  by two. Precisely,
by general properties of double covers, there is an invertible
sheaf $\delta/2$ on $\cur_2$ satisfying:
\[
\sigma_*\O_{\cur\times\cur} = \O_{\cur_2}\oplus
\O_{\cur_2}(-\delta/2)
\]
and
\begin{equation}\label{E:sim1}
\sigma^*(\delta/2) = D
\end{equation}
In particular we have:
\begin{equation}\label{E:sim1.5}
(\delta/2\cdot \delta/2) = 1-\gamma, \qquad \O_\delta(-\delta/2) =
\omega_\delta
\end{equation}
   The \emph{numerical class}   of   the
irreducible curve in $\cur_2$:
\[
x+ \cur = \{x+y: y\in \cur\}\subset \cur_2
\]
  is independent of  $x \in \cur$ and  will be denoted  by $X$.
If $L=\O(\sum n_ix_i)$  is an invertible sheaf on $\cur$ we will
denote by $LX$ the invertible sheaf $\O_{\cur_2}(\sum
n_i(x_i+\cur))$ on $\cur_2$.  If deg$(L)=n$ then
\[
LX \equiv_{num} nX
\]
It is immediate to check that
\begin{equation}\label{E:sim2}
\sigma^*(LX) = L\boxtimes L
\end{equation}
and
\begin{equation}\label{E:sim3}
X\cdot \delta/2 = 1 = X^2,  \qquad
\end{equation}

\begin{lemma}\label{L:sim1}
\[
\omega_{\cur_2} = \omega_\cur X - \delta/2 
\]
\end{lemma}

\proof It is a special case of a well-known formula of Mattuck
\cite{aM65}. Let's give a direct proof.  Since we have
\[
\omega_{\cur\times \cur} = \omega_\cur\boxtimes\omega_\cur
\]
by the Hurwitz formula for the double cover $\sigma$
(\cite{BPVdV}, (19) p. 41) we obtain:
 \[
  \omega_\cur\boxtimes\omega_\cur(-D) = \sigma^*\omega_{\cur_2}
  \]
On the other hand, by (\ref{E:sim1}) and (\ref{E:sim2})  we have:
\[
 \sigma^*\left[\omega_\cur X - {\delta\over 2}\right] =  \omega_\cur\boxtimes\omega_\cur(-D)
 \]
 In other words:
 \[
 \sigma^*\omega_{\cur_2} =  \sigma^*\left[\omega_\cur X - {\delta\over 2}\right]
 \]
 Since $\sigma^*: {\rm Pic}(\cur_2) \to {\rm Pic}(\cur\times \cur)$ is injective it follows that:
 \[
 \omega_{\cur_2} = \omega_\cur X - {\delta\over 2}
 \]
 \qed

 \begin{remark}\label{R:sim1}\rm
 Comparing (\ref{E:sim1.5}) with the adjunction formula we obtain:
 \[
 \O_\delta(-\delta/2) = \O_\delta(\omega_\cur X+ \delta/2)
 \]
 and therefore $\O_\delta(\omega_\cur X+\delta) = \O_\delta$. This implies that the effective divisor
 $\omega_\cur X+\delta$ on $\cur_2$ is not ample.
 \end{remark}

 \begin{lemma}\label{L:sim0}
Let $L \in {\rm Pic}(\cur)$. Then
\[
\begin{array}{llll}
H^0(\cur_2,LX) = &S^2H^0(\cur,L) \\ \\
H^1(\cur_2,LX) \cong &H^0(\cur,L)\otimes H^1(\cur,L) \\ \\
H^2(\cur_2,LX) = &\bigwedge^2 H^1(\cur,L)
\end{array}
\]
and
\[
\begin{array}{llll}
H^0(\cur_2,LX-\delta/2) =& \bigwedge^2H^0(\cur,L) \\ \\
H^1(\cur_2,LX-\delta/2) \cong &H^0(\cur,L)\otimes H^1(\cur,L) \\ \\
H^2(\cur_2,LX-\delta/2) = &S^2 H^1(\cur,L)
\end{array}
\]
\end{lemma}

\proof The involution $\iota: \cur\times\cur \to \cur\times\cur$
acts on
\[
H^0(\cur\times\cur, L\boxtimes L)=H^0(\cur,L) \otimes H^0(\cur,L)
= S^2H^0(\cur,L)\oplus \bigwedge^2H^0(\cur,L)
\]
by interchanging the factors in the tensor product.   On the other
hand, by the projection formula we have
\[
H^0(\cur\times\cur, L\boxtimes L) = H^0(\cur\times\cur,
\sigma^*(LX)) = H^0(\cur_2,LX)\oplus H^0(\cur_2,LX-\delta/2)
 \]
  and the two summands are respectively the $\iota$-invariant and $\iota$-antinvariant part.
  This proves the equalities for $H^0$'s. The equalities
  for $H^2$'s are proved similarly using Serre duality and the expression of $\omega_{\cur_2}$.

  \noindent
  The isomorphisms for  $H^1$'s follow from the Kunneth formula and  direct computation.  \qed

\begin{lemma}\label{L:sim2}
There are natural identifications:
\begin{description}
\item{(i)} \[ H^1(\cur_2,\O_{\cur_2}) = H^1(\cur,\O_\cur)
\]
\item{(ii)}  \[ H^2(\cur_2,\O_{\cur_2}) = \bigwedge^2
H^1(\cur,\O_\cur) = H^0(\cur_2,\omega_{\cur_2})^\vee
\]
\item{(iii)}
\[H^1(\cur_2,\omega_{\cur_2}) = H^0(\cur,\omega_\cur)
\]
 \end{description}
 \end{lemma}

 \proof
  It is a special case of Lemma \ref{L:sim0}.  \qed

  \begin{remark}\label{R:sim2}\rm
  After making the natural identification $\delta=\Gamma$, from Lemma \ref{L:sim0} and
  from (\ref{E:sim1.5}) and (\ref{E:sim3})  it follows that the restriction map:
  \[
  \xymatrix{
  H^0(\cur_2,LX)\ar[r] & H^0(\delta,LX\otimes \O_\delta)}
  \]
  is identified  with the natural restriction:
  \[
  \xymatrix{S^2H^0(\cur,L) \ar[r] & H^0(\cur,L^2)}
  \]
  and therefore
  \[
  H^0(\cur_2,LX-\delta) = \ker\left[\xymatrix{S^2H^0(\cur,L) \ar[r] & H^0(\cur,L^2)}\right]
  \]
  Similarly the map:
  \[
   \xymatrix{
  H^0(\cur_2,LX-\delta/2)\ar[r] & H^0(\delta,(LX-\delta/2)\otimes \O_\delta)}
  \]
  is identified with the Wahl map:
  \[
  \xymatrix{
  w_L:\bigwedge^2 H^0(\cur,L) \ar[r] & H^0(\cur,\omega_\cur\otimes L^2)}
  \]
and $H^0(\cur_2,LX-3\delta/2)=\ker(w_L)$.
  \end{remark}

Consider the universal divisor
$$
\Delta := \{(x,x+y): x,y \in \cur\} \subset \cur\times \cur_2
$$
and   the diagram:
\[
\xymatrix{\Delta \ar@{^(->}[r]&\cur \times \cur_2 \ar[dr]^{q} \ar[dl]_p \\
\cur && \cur_2}
\]
We have an isomorphism:
\[
\xymatrix{\epsilon: \cur\times \cur \ar[r] & \Delta \\
(x,y) \ar@{|->}[r] & (x,x+y)}
\]
and we can identify the map $\sigma$ with the composition
\[
\xymatrix{\cur\times\cur \ar[r]^\epsilon & \Delta \ar[r]^q &
\cur_2}
\]

 Given    a line bundle $L$  on $\cur$
 we can consider the  locally free sheaf of rank two on $\cur_2$:
 \[
 E_L := q_*(p^*L_{|\Delta})
 \]
which is called the \emph{secant bundle} of $L$.  Since
\[
\xymatrix{ p\cdot\epsilon = \pi_1: \cur\times\cur \ar[r] & \cur}
\]
we may as well write:
\[
E_L = \sigma_*(\pi_1^*L)
\]

\noindent
 Let  $(L,V)$ be a $g^1_n$  on  $\cur$.
  On $\cur_2$ consider the composition
  \[
  \xymatrix{
  \phi: V\otimes \O_{\cur_2}\ar@{^(->}[r]&q_*p^*L \ar[r] & E_L}
  \]
  where the second homomorphism is obtained by pushing down
  the restriction homomorphism on $\cur\times\cur_2$:
  \[
  \xymatrix{
   p^*L \ar[r] & p^*L_{|\Delta}}
  \]
  $\phi$ is a homomorphism of locally free sheaves of rank
   two. We define an effective divisor $s_{L,V}$ on $\cur_2$ by:
  \[
  s_{L,V} := D_1(\phi) = D_0(\bigwedge^2\phi)
  \]
where
\[
\xymatrix{ \bigwedge^2\phi: \bigwedge^2V\otimes\O_{\cur_2} \ar[r]
& \bigwedge^2E_L}
\]
 Set-theoretically we have:
\[
s_{L,V} = \left\{x+y\in \cur_2: \dim[V(-x-y)] \ge 1\right\}
\]

\begin{lemma}\label{L:sLV1}
Let  $(L,V)$ be a $g^1_n$  on  $\cur$.
\begin{itemize}
\item[(i)] We have
\[
\O(s_{L,V}) = \bigwedge^2V^\vee \otimes \bigwedge^2E_L \cong LX-
\delta/2
\]

\item[(ii)]
\[
\begin{array}{lll}
H^0(\cur_2,\O_{\cur_2}(s_{L,V})) = & \bigwedge^2V^\vee \otimes\bigwedge^2 H^0(\cur,L) \\ \\
H^1(\cur_2, \O_{\cur_2}(s_{L,V})) = &\bigwedge^2V^\vee \otimes H^0(\cur,L)\otimes H^1(\cur,L) \\ \\
 H^2(\cur_2, \O_{\cur_2}(s_{L,V})) = &\bigwedge^2V^\vee \otimes S^2H^1(\cur,L)
 \end{array}
 \]
\item[(iii)]  The arithmetic genus $g(s_{L,V})$ of $s_{L,V}$
satisfies:
\[
2g(s_{L,V}) - 2 = (n-2)(n+2\gamma-3)-2
\]
and
\[s_{L,V}\cdot s_{L,V} = (n-1)^2- \gamma\]
\end{itemize}
\end{lemma}

\proof (i) See \cite{wF84}, Ex. 14.4.17, p. 263. (ii) follows from
Lemma \ref{L:sim0}. (iii) is left to the reader.
 \qed

 Let  $(L,V)$ be a $g^1_n$  on  $\cur$, and let:
 \[
  t_{L,V} := \{(x,y): V(-x-y) \ne 0\} \subset \cur\times\cur
  \]
   Consider the morphisms:
\begin{equation}\label{E:square1}
\xymatrix{
&t_{L,V}\ar[dl]_-{\pi_1}\ar[dr]^-\sigma  \\
\cur&&s_{L,V}}
\end{equation}
  where $\pi_1$ is induced by the first projection.
  Assume that $(L,V)$ is base-point free. Then the above diagram can be completed as:
  \begin{equation}\label{E:square2}
\xymatrix{
&t_{L,V}\ar[dl]_-{\pi_1}\ar[dr]^-\sigma \\
\cur\ar[dr]_-{\varphi_{L,V}}&&s_{L,V} \ar[dl]^-f \\ &\P(V^\vee)}
\end{equation}
where  $f$ is the obvious morphism.  The degrees of the morphisms
are:
\[
\deg(\pi_1)=n-1, \ \deg(\varphi_{L,V})=n, \ \deg(\sigma)=2, \
\deg(f)= {n\choose 2}
\]

 \begin{lemma}\label{L:sLV2}
 Let  $(L,V)$ be a $g^1_n$  on  $\cur$. Then:

 \begin{itemize}
 \item[(i)]   $t_{L,V}$ and $s_{L,V}$ are both connected.

\item[(ii)] If  $(L,V)$ is  base-point free and with simple
ramification then $t_{L,V}$ and $s_{L,V}$ are both irreducible and
nonsingular.

\item[(iii)] Assume that $\cur$ is  very general and that  $n \leq
\gamma -1$. If $s_{L,V}$ is reducible then the fixed  divisor $B$
of $(L,V)$ is positive
   and $s_{L,V}=BX \cup s_{L(-B),V(-B)}$.
\end{itemize}
\end{lemma}

\proof (i) Let's prove that $s_{L,V}$ is connected. We may assume
that $n \ge 2$. Since $s_{L,V}$ is effective we have
$H^0(\O_{\cur_2}(-s_{L,V}))=0$. Therefore from the exact sequence
on $\cur_2$:
\[
\xymatrix{ 0 \ar[r] & \O(-s_{L,V}) \ar[r]& \O_{\cur_2}\ar[r] &
\O_{s_{L,V}}\ar[r]&0}
\]
 we see that  it suffices to prove that
$H^1(\O_{\cur_2}(-s_{L,V}))=0$. Consider the exact sequence:
\[
\xymatrix{ 0 \ar[r]& \O(-LX-{\delta\over 2}) \ar[r]& \O(-s_{L,V})
\ar[r]& \O_\delta(-s_{L,V}) \ar[r]&0}
\]
We have $H^1(\O(-LX-{\delta\over 2}) = H^0(-L)\otimes H^1(-L)=0$
(Lemma \ref{L:sim0}). Thus it suffices to prove that
$\xymatrix{\partial:H^1(
\O_\delta(-s_{L,V}))\ar[r]&H^2(\O(-LX-{\delta\over 2})}$ is
injective, or equivalently  that  $\partial^\vee$ is surjective.
But $\O_\delta(-s_{L,V})=\O_\cur(-\omega_\cur-2L)$ and
$H^2(\O(-LX-{\delta\over 2})=S^2H^1(\cur,-L)$ (Lemma
\ref{L:sim0}).  Therefore
\[
\xymatrix{
\partial^\vee: S^2H^0(\omega_\cur+L) \ar[r] & H^0(2(\omega_\cur+L))}
\]
This map is the natural multiplication, and it is surjective
because $\mathrm{deg}(\omega_\cur+L) \ge 2g+1$ (\cite{dM70}). This
concludes the proof of the connectedness of  $s_{L,V}$.

\noindent Since $\sigma$ is ramified over the points $2x \in
s_{L,V}$, where $x \in \cur$ is a ramification point of
$\varphi_{L,V}$, we deduce that $t_{L,V}$ is connected as well.

(ii)   Assume that $(L,V)$ is base-point free. Let
Sing$(t_{L,V})\subset t_{L,V}$ and Sing$(s_{L,V})\subset s_{L,V}$
be the singular loci. Since $\cur$ and $\P(V^\vee)$ are
nonsingular we have:
 \[
 \begin{array}{ll}
 {\rm Sing}(s_{L,V})\subset {\rm Ram}(f) &= \{x+y\in s_{L,V}: \textrm{either $x$ or $y$ is a ramif. pt of $(L,V)$}\} \\ \\
 {\rm Sing}(t_{L,V}) \subset {\rm Ram}(\pi_1)&= \{(x,y): V(-x-2y)\ne 0 \}
 \end{array}
 \]
Assume now that $(L,V)$ has only simple ramification and that
$(x,y)\in {\rm Ram}(\pi_1)$.
 Then $x\ne y$ because otherwise $V(-3x)\ne 0$. Moreover $(y,x) \notin {\rm Ram}(\pi_1)$ because otherwise
$V(-2x-2y)\ne 0$. Therefore $(y,x)\notin {\rm Sing}(t_{L,V})$. It
follows that $\sigma(x,y)=\sigma(y,x) \notin {\rm Sing}(s_{L,V})$
because $\sigma$ is etale over $x+y$. But then also $(x,y)\notin
{\rm Sing}(t_{L,V})$. The conclusion is that   $t_{L,V}$ and
$s_{L,V}$ are both nonsingular. Since they are connected they are
irreducible as well.

(iii)
 If $\cur$ is very general, from [ACGH] Lemma p. 359, it follows that NS$(\cur_2)=<x,\delta/2>$, so that
 in particular the numerical class of every effective
 curve is a combination of $x$ and $\delta/2$ with integer coefficients.
 Let us assume that $s_{(L,V)}=A \cup B$, where $A$ and $B$ are effective curves on $\cur_2$.

 We write $A =a_1x+a_2\delta/2$, $B =b_1x+b_2\delta/2$ for their numerical classes.
 We have $a_2+b_2=-1$ and we will assume that $a_2 \geq b_2$. \par \noindent
 Suppose that $a_2 \geq 1$, so that $b_2 \leq -2$. Since $\delta$ is irreducible and is not contained in any
 component of $s_{(L,V)}$   we have $$\delta \cdot A \geq 0$$
 so that $$a_1 \geq a_2(\gamma - 1) \geq n.$$
 It then follows that $b_1=n-a_1 \leq 0$. But then $x \cdot B=b_1+b_2 \leq -2$, and $x$ is the class
 of an ample divisor. \par \noindent

 It then follows that $a_2 \leq 0$, so that in fact $a_2=0$ because $a_2+b_2=-1$ and $a_2 \geq b_2$.
 Then $A=a_1x$ numerically and in particular $(L,V)$ has a base locus of degree $a_1$. If $B$ is reducible
 we repeat the argument until we get down to a base point free pencil.
 \qed

\begin{remarks}

\end{remarks}

\indent (a)  If $(L,V)$ does not have simple ramification then
$t_{L,V}$ and $s_{L,V}$ can be singular. Assume for example that
the pencil contains a divisor of the form
$2x+2y+z_1+\cdots+z_{n-4}$, with $x\ne y$. Then $V(-x-2y)\ne 0 \ne
V(-2x-y)$. This implies that $(x,y)\in t_{L,V}$ is in
$\mathrm{Ram}(\pi_1) \cap \mathrm{Ram}(\pi_2)$, where
$\pi_1,\pi_2: \cur\times\cur \longrightarrow \cur$ are the
projections. Therefore $(x,y)\in \mathrm{Sing}(t_{L,V})$, and
$\sigma(x,y)\in \mathrm{Sing}(s_{L,V})$.

(b) Part (iii) of Lemma \ref{L:sLV2} will not be needed in the
rest of the paper.

(c)  The curves $t_{L,V}$ are also considered in the recent
preprint \cite{GK11}, where they are called \emph{trace curves.}

\section{The schemes of neutral linear series}\label{S:NLS}

We shall adopt the standard notation $G^r_n(\cur)$ to denote the
scheme of $g^r_n$'s on the curve $\cur$. For its definition and
main properties we refer to  \cite{ACGH}.

\begin{definition}\label{D:nls1}
Let $(L,V)$ be a $g^r_n$ on $\cur$ and $x+y \in \cur_2$. We say
that $(L,V)$ is neutral w.r. to $x+y$ if $\mathrm{dim}(V(-x-y))
\ge r$.
\end{definition}

  Consider the
following diagram where all the arrows are projections:
\begin{equation}\label{E:big}
\xymatrix{ \cur \times \cur_2 \ar[d]_q& \cur \times \cur_2 \times
G^r_n \ar[r]^-{\pi_{13}}\ar[d]^-{\pi_{23}} \ar[l]_-{\pi_{12}} &
 \cur \times G^r_n
\ar[d]^-u\\
 \cur_2&\cur_2 \times G^r_n \ar[l]_-{\pi_1} \ar[r]^-{\pi_2}&G^r_n \\
}
\end{equation}
and where   $G^r_n = G^r_n(\cur)$.   Let $\cP$ be the pullback on
$\cur \times G^r_n$
 of a Poincar\'e line bundle on $\cur \times \Pic^n(\cur)$ and let $E^r_n \subset u_*\cP$
 be the tautological locally free   subsheaf of rank $r+1$ on $G^r_n$.

   We   consider  on $\cur \times \cur_2 \times G^r_n$ the natural restriction homomorphism:
\[
\xymatrix{ \pi_{13}^*\cP \ar[r] & \pi_{13}^*\cP \otimes
\pi_{12}^*\O_\Delta}
\]
Then   we push it down to $\cur_2 \times G^r_n$ and we compose
with the inclusion
\[
\pi_2^*E^r_n \subset \pi_2^*u_*\cP=\pi_{23*}\pi_{13}^*\cP \]
 We
obtain a homomorphism of locally free   sheaves on $\cur_2 \times
G^r_n$:
\[
\xymatrix{ \Phi: \pi_2^*E^r_n \ar[r] & \pi_{23*}[\pi_{13}^*\cP
\otimes \pi_{12}^*\O_\Delta]}
\]
of ranks $r+1$ and two respectively. Consider the closed
subscheme:
\[
N^r_n(\cur) := D_1(\Phi) \subset \cur_2 \times G^r_n
\]
and the projections:
\begin{equation}\label{E:proj1}
\xymatrix{
N^r_n(\cur)\ar[d]_-{p_1}\ar[dr]^-{p_2} \ar@{^(->}[r] & \cur_2 \times G^r_n \\
\cur_2&G^r_n}
\end{equation}
 For each $x+y \in \cur_2$ we have
\[
p_1^{-1}(x+y)= D_1(\Phi(x+y)) \subset G^r_n
\]
   where
 \[
 \xymatrix{
 \Phi(x+y): E^r_n  \ar@{^(->}[r]&u_*\cP \ar[r] & \cP\otimes \O_{x+y}}
 \]
  We define
  \[
  N^r_n(x+y):= p_1^{-1}(x+y)
  \]
  Set-theoretically:
\[
N^r_n(x+y) = \{(L,V): \dim[V(-x-y)] \ge r\}
\]
Therefore we call  $N^r_n(x+y)$ \emph{the scheme of neutral linear
$g^r_n$'s with respect to $x+y$}. If $r=1$ then  for every $(L,V)
\in G^1_n$ we have an isomorphism:
\[
p_2^{-1}(L,V) \cong s_{L,V}
\]
In particular we see that $N^1_n(\cur)\ne \cur_2 \times G^1_n$ and
therefore it is a proper divisor in $\cur_2 \times G^1_n$.

\section{The case of pencils}\label{S:NLPII}

  \emph{ In this section we  assume that $\cur$ is Petri w.r. to pencils and that}
 \begin{equation}\label{E:rho1}
 \gamma \ge\rho(\gamma,1,n) := 2n-2-\gamma \ge 1
 \end{equation}

\begin{lemma}\label{L:rho1}
Under the assumptions above   $G^1_n=G^1_n(\cur)\ne \emptyset$ is
nonsingular connected of   dimension $\rho(\gamma,1,n)$, and
consists generically of complete and base point free $g^1_n$'s.
\end{lemma}

\proof  For the first assertion see \cite{ACGH}, Prop. 4.1 p. 187.
If a component of $G^1_n$ consists generically of $g^1_n$'s with
base points, then $G^1_{n-1}$ has a component of dimension $\ge
\rho(\gamma,1,n)=\rho(\gamma,1,n-1)+2$, contradicting the
hypothesis that $\cur$ is Petri w.r. to pencils. A similar
argument holds in case a component of $G^1_n$ consists of
incomplete pencils. \qed

 We have a natural inclusion:
\[
\xymatrix{ \cur_2\times G^1_{n-2}\ar[r]^-\iota & \cur_2\times
G^1_n
\\ (x+y,(A,U)) \ar@{|->}[r]& (x+y,(A(x+y),U(x+y))}
\]

\begin{lemma}\label{L:D0}
\[\iota(\cur_2\times G^1_{n-2}) = D_0(\Phi) \subset
\mathrm{Sing}(N^1_n(\cur))
\]
\end{lemma}

\proof $\iota(\cur_2\times G^1_{n-2})$ is supported on the points
$(x+y,(L,V))\in \cur_2\times G^1_n$ such that $x+y$ is contained
in the fixed divisor of $(L,V)$, and these are exactly the points
where $\Phi$ vanishes. The inclusion is obvious. \qed

Note that  $\cur_2\times G^1_{n-2}=\emptyset$ if $\rho(\gamma,1,n)
\le 3$, and it is nonsingular of pure dimension
\[
\dim(\cur_2\times
G^1_{n-2})=\rho(\gamma,1,n-2)+2=\rho(\gamma,1,n)-2 \]
 if
$\rho(\gamma,1,n)\ge 4$.

   \begin{lemma}\label{L:kerM}
   Let $(L,V) \in N^1_n(x+y)$ for some $x+y\in \cur_2$, and consider the multiplication map:
   \[
   \xymatrix{
   M: V \otimes H^0(\omega_\cur L^{-1}(x+y))\ar[r] & H^0(\omega_\cur(x+y))}
   \]
   Let $B$ be the fixed divisor of $(L,V)$. Then
   \begin{description}
   \item[(a)]
   $ \mathrm{ker}(M)\cong H^0(\omega_\cur L^{-2}(x+y+B)) \cong H^1(L^2(-x-y-B))^\vee
   $
 \item[(b)]  $\mathrm{dim}[\mathrm{ker}(M)] \le 1$.

  \item[(c)]
   If $x < B$   then $M$ is injective.
   \end{description}
      \end{lemma}

   \proof
   (a) is  a direct consequence of (\ref{E:LM3}) and (b) is left to the reader.

   \noindent
   (c) From Proposition \ref{P:kermu}(i) it follows that $H^1(L^2(-2B))=0$, and   then also
   $H^1(L^2(-2B-y))=0$, by Prop. \ref{P:kermu}(ii).
     But, since $x < B$, we have a surjection
     \[
     \xymatrix{
     H^1(L^2(-2B-y))\ar@{->>}[r] & H^1(L^2(-x-y-B))\cong \mathrm{ker}(M)^\vee}
     \]
      and therefore $\mathrm{ker}(M)=0$.  \qed

   \begin{prop}\label{P:TN}
   Let $(x+y,(L,V)) \in N^1_n(\Gamma)$.
   Then $(L,V)$ is a nonsingular point
   of $N^1_n(x+y)$ if and only if $(x+y,(L,V))\notin \iota(\cur_2\times G^1_{n-2})$ and $M$ is injective.
   \end{prop}

   \proof
    For every $x+y \in \cur_2$ the points of $N^1_n(x+y) \cap
\iota(\cur_2\times G^1_{n-2})$ are singular for $N^1_n(x+y)$, by
Lemma \ref{L:D0}. Therefore we can assume that $(x+y,(L,V)) \in
N^1_n(\Gamma)\setminus \iota(\cur_2\times G^1_{n-2})$.
  We need to describe the tangent space
$T_{(L,V)}N^1_n(x+y) \subset T_{(L,V)}G^1_n$. Note that
   \begin{equation}\label{E:dimT}
   \rho-1 \le \dim\left[T_{(L,V)}N^1_n(x+y)\right] \le \rho
   \end{equation}
   where $\rho:=\rho(\gamma,1,n)$
   and the first equality holds if and only if  $N^1_n(x+y)$
   is nonsingular of dimension $\rho-1$ at $(L,V)$.
    Denote by
   \[
  \xymatrix{
   \pi: G^1_n \ar[r] & \Pic^n(\cur)}
   \]
   the Abel-Jacobi map.
   Consider the well-known exact sequence (\cite{ACGH}, p. 187):
\begin{equation}\label{E:TG}
\xymatrix{
&&T_{(L,V)}N^1_n(x+y)\ar@{^(->}[d] \\
0\ar[r] & \Hom(V,H^0(L)/V) \ar[r]& T_{(L,V)}G^1_n\ar[r]^-{d\pi}
&H^1(\O_\cur) \ar[r]^-{\mu_0^\vee} & V^\vee\otimes H^1(L)}
\end{equation}
From this sequence and from (\ref{E:dimT}) we see that $(L,V)$ is
a nonsingular point of $N^1_n(x+y)$ if and only if either $
\Hom(V,H^0(L)/V) \not\subset T_{(L,V)}N^1_n(x+y)$ or
$d\pi\left(T_{(L,V)}N^1_n(x+y)\right) \ne
d\pi\left(T_{(L,V)}G^1_n\right)$.

\emph{Assume first that  the complete linear series $(L,H^0(L))$
is not neutral w.r. to $x+y$;}  then in particular $V \ne H^0(L)$
and $ \Hom(V,H^0(L)/V) \not\subset T_{(L,V)}N^1_n(x+y)$. It
follows that $N^1_n(x+y)$ is nonsingular at $(L,V)$. On the other
hand we have $H^0(\omega_\cur L^{-1}(x+y)) \cong H^0(\omega_\cur
L^{-1})$ and from the commutative diagram:
\[
\xymatrix{
V\otimes H^0(\omega_\cur L^{-1}) \ar[r]^-{\mu_0} \ar[d]& H^0(\omega_\cur) \ar[d] \\
V\otimes H^0(\omega_\cur L^{-1}(x+y))\ar[r]^-M &
H^0(\omega_\cur(x+y))}
\]
we see that $M$ is injective. Therefore the Proposition is true in
this case.

\emph{Assume    now    that the complete series $(L,H^0(L))$ is
neutral w.r. to $x+y$.}
 In this case
 \[
  \Hom(V,H^0(L)/V) \subset T_{(L,V)}N^1_n(x+y)
 \]
 and
 therefore $(L,V)$ is a nonsingular point of $N^1_n(x+y)$ if and only if
\[
d\pi\left(T_{(L,V)}N^1_n(x+y)\right) \ne
d\pi\left(T_{(L,V)}G^1_n\right)=  \ker(\mu_0^\vee)
\]
 We have dim$[V(-x-y)]=1$ because $(x+y,(L,V)) \notin
   D_0(\Phi)$. Let $\mu =
\mu_0(L(-x-y),V(-x-y))$, so that
\[
\xymatrix{ \mu^\vee: H^1(\O_\cur)\ar[r] & V(-x-y)^\vee\otimes
H^1(L(-x-y))}
\]
 and $\ker(\mu^\vee)$ consists of the tangent directions along which
 $L(-x-y)$ deforms carrying along the   section of $V(-x-y)$.

Then
   \begin{equation}\label{E:petri0}
 d\pi\left[T_{(L,V)}N^1_n(x+y)\right] = \ker(\mu_0^\vee)\cap \ker(\mu^\vee)
 \end{equation}
 Consider the following   commutative diagram:
 \[
 \xymatrix{
 &V^\vee\otimes H^1(L) \ar[dr] \\
 H^1(\O_\cur)\ar[ur]^-{\mu_0^\vee}\ar[r]^-\nu \ar[dr]_-{\mu^\vee} & P \ar[u]\ar[d] & V(-x-y)^\vee\otimes H^1(L) \\
 &V(-x-y)^\vee\otimes H^1(L(-x-y) \ar[ur]}
 \]
  where
\[
P := [V^\vee \otimes H^1(L)] \times_{V(-x-y)^\vee\otimes H^1(L)}
[V(-x-y)^\vee\otimes H^1(L(-x-y))]
\]
and where $\nu$ is the map induced by the universal property of
$P$.  By construction
\begin{equation}\label{E:petri1}
  \ker(\nu)= \ker(\mu_0^\vee)\cap\mathrm{ker}(\mu^\vee)
\end{equation}
   Observing that $h^1(L(-x-y))=h^1(L)+1$, we deduce that dim$(P)=2h^1(L)+1$.
   Therefore, recalling (\ref{E:petri0}) and \ref{E:petri1}), we
   see that \emph{$(L,V)$ is a nonsingular point of
$N^1_n(x+y)$ if and only if $\nu$ is surjective.}

   We   have another commutative diagram:
  \[
  \xymatrix{
 &V^\vee\otimes H^1(L) \ar[dr] \\
 V^\vee\otimes H^1(L(-x-y))\ar[ur]\ar[r]^-n \ar[dr]  & P \ar[u]\ar[d] & V(-x-y)^\vee\otimes H^1(L) \\
 &V(-x-y)^\vee\otimes H^1(L(-x-y)) \ar[ur]}
 \]
  where $n$ is induced by the two surjective maps $e:H^1(L(-x-y)) \to H^1(L)$ and
  $p:V^\vee \to V(-x-y)^\vee$. It follows from linear algebra that $n$  is surjective
  and that
  $\ker(n) = \ker(e)\otimes \ker(p)$ is one-dimensional.
  Therefore there is induced the following commutative and exact diagram:
     \begin{equation}\label{E:petri2}
  \xymatrix{
  H^1(\O_\cur(-x-y)) \ar[d]_-\epsilon \ar[r]^-{M^\vee}&V^\vee\otimes H^1(L(-x-y)) \ar[r]\ar[d]^-n &
  \mathrm{coker}(M^\vee)\ar[r] \ar[d]&0 \\
  H^1(\O_\cur) \ar[r]^\nu& P \ar[r] &\mathrm{coker}(\nu) \ar[r] & 0}
  \end{equation}
 where the vertical map on the right is surjective. Therefore, if $M^\vee$ is surjective
 then $\nu$ is surjective. Assume that $M^\vee$ is not surjective.
 Then   neither $x$ nor $y$ is a  base points of $(L,V)$
 (Lemma  \ref{L:kerM}(c)) and    Im$(M) \not\subset H^0(\omega_\cur)$:   it follows that
$\ker(M^\vee) \cap \ker(\epsilon)=0$, and this implies that
$(0)\ne \mathrm{coker}(M^\vee)\cong \mathrm{coker}(\nu)$. In
conclusion  $M^\vee$ is surjective if and only if $\nu$ is
surjective
 if and only if  $(L,V)$ is a nonsingular point of
$N^1_n(x+y)$, and this proves
 the Proposition.\qed

 \begin{corollary}\label{C:petri3}
 Assume that $(x+y,(L,V))\in N^1_n(\cur)\setminus \iota(\cur_2\times G^1_{n-2})$.
 Let $B$ be the fixed divisor
  of $(L,V)$ and $b:=\mathrm{deg}(B)$.
   The following conditions are equivalent:
 \begin{itemize}

 \item[(i)] $h^1(L^2(-x-y-B)) =1$.

 \item[(ii)] There exists a unique effective  divisor
 $z_1+\cdots + z_{2\gamma-2n+b}\in |\omega_\cur L^{-2}(x+y+B)|$.

 \item[(iii)]
   $(L,V)$ is a singular point of $N^1_n(x+y)$.
   \end{itemize}
   In particular, if $x<B$ then $(L,V)$ is a nonsingular point of $N^1_n(x+y)$.
    \end{corollary}

 \proof
  The equivalence of (i) and (ii) is obvious.   By Proposition \ref{P:TN}
 condition (iii) is satisfied if and only if $M$ is not injective, and this is equivalent to
 $\dim(\mathrm{ker}(M))=h^1(L^2(-x-y-B)) =1$, by Lemma
 \ref{L:kerM}. The last assertion is Lemma \ref{L:kerM}(c).
    \qed

Let's consider the diagram (\ref{E:proj1}) in the case $r=1$:
\begin{equation}\label{E:proj2}
\xymatrix{
N^1_n(\cur)\ar[d]_-{p_1}\ar[dr]^-{p_2} \ar@{^(->}[r] & \cur_2 \times G^1_n \\
\cur_2&G^1_n}
\end{equation}
The projection $p_2$ defines a flat family of curves in $\cur_2$,
namely the family of all curves $s_{L,V}$. Thus, by functoriality,
we have a morphism:
\[
\xymatrix{ \zeta: G^1_n \ar[r] & \mathrm{Hilb}^{\cur_2}}
\]

\begin{theorem}\label{T:hilb}
    The morphism $\zeta$ is an isomorphism of $G^1_n$ onto an
    irreducible component of $\mathrm{Hilb}^{\cur_2}$.
\end{theorem}

\proof $\zeta$ is injective: in fact $(L,V)$ can be reconstructed
from $s_{L,V}$ as the pencil generated by the two divisors
$(a+\cur)\cap s_{L,V}$ and $(b+\cur)\cap s_{L,V}$ for distinct
general $a,b \in \cur$. Therefore, since $G^1_n$ is nonsingular,
it will suffice to prove that $d_{(L,V)}\zeta$ is   injective for
all $(L,V)$ and that it is an isomorphism in case   $(L,V)$ is
complete.

Recalling Lemmas \ref{L:sim0} and  \ref{L:sLV1} we see that the
cohomology sequence associated to the exact sequence
\[
\xymatrix{ 0 \ar[r]&\O_{\cur_2}\ar[r] & \O_{\cur_2}(s_{L,V})\ar[r]
& N_{s_{L,V}} \ar[r] &0 }
\]
is:
\begin{equation}\label{E:normal}
    \xymatrix{
    0 \ar[r] &
    \frac{\bigwedge^2V^\vee\otimes\bigwedge^2H^0(\cur,L)}{H^0(\cur,\O_{\cur})} \ar[r]&
    H^0(s_{L,V},N_{s_{L,V}})\ar[r]&
    H^1(\cur,\O_\cur)\ar[r]&\bigwedge^2V^\vee\otimes H^0(\cur,L)\otimes H^1(\cur,L)
    }
\end{equation}
 Comparing
with the exact sequence (\ref{E:TG}) we deduce a commutative
diagram:
\[
\xymatrix{ 0\ar[r] & \Hom(V,H^0(L)/V) \ar[r]\ar@{^(->}[d]&
T_{(L,V)}G^1_n\ar[r]^-{d\pi}
\ar[d]^-{d_{(L,V)}\zeta}&H^1(\O_\cur)\ar@{=}[d]
\ar[r]^-{\mu_0(L,V)^\vee} &
V^\vee\otimes H^1(L) \\
0 \ar[r] &
    \frac{\bigwedge^2H^0(\cur,L)}{H^0(\cur,\O_{\cur})} \ar[r]& H^0(s_{L,V},N_{s_{L,V}})\ar[r]&
    H^1(\O_\cur)\ar[r]&\bigwedge^2V^\vee\otimes H^0(\cur,L)\otimes H^1(\cur,L)
    }
\]
where the left vertical arrow is the differential at $V\subset
H^0(L)$ of the Pl\"ucker embedding of $\mathrm{Grass}(2,H^0(L))$.
Clearly $d_{(L,V)}\zeta$ is injective and it is an isomorphism
when $V=H^0(L)$. \qed

  \section{The singularities of $p_1$}\label{S:singp1}

 We keep the same assumptions as in \S \ref{S:NLPII}, namely that $\cur$ is Petri w.r. to pencils
 and that (\ref{E:rho1}) holds.
 Let
\[
\I := \O(-N^1_n(\cur)) \subset \O_{\cur_2\times G^1_n}
\]
be the ideal sheaf of $N^1_n(\cur)$ inside $\cur_2\times G^1_n$.
Consider the natural homomorphism of locally free sheaves on
$N^1_n(\cur)$:
\[
\xymatrix{ \partial:\I/\I^2 \ar[r] & p_2^*\Omega^1_{G^1_n} =
\Omega^1_{\cur_2\times G^1_n/\cur_2|N^1_n}}
\]
We will denote by $D_0(\partial)$ the vanishing scheme of
$\partial$.

\begin{prop}\label{P:NP1}
Assume that $\cur$ is Petri w.r. to pencils.  Then
\begin{description}
\item[(i)] Supp$(D_0(\partial))$ coincides with the relative
singular locus   of $\xymatrix{p_1:N^1_n(\Gamma) \ar[r]& \cur_2}$.
In particular $\iota(\cur_2\times
 G^1_{n-2})\subset \mathrm{Supp}(D_0(\partial))$.
\item[(ii)] Let $U \subset G^1_n$ be the open set of $g^1_n$'s
which are base point free and complete. Then $D_0(\partial)\cap
p_2^{-1}(U) \ne \emptyset$ and every irreducible component of $F:=
\overline{D_0(\partial)\cap p_2^{-1}(U)}$ has dimension $\ge 1$.
\end{description}
 \end{prop}

\proof
 (i)  Since we have
\[
\Omega^1_{\cur_2\times G^1_n/\cur_2} = \pi^*_2\Omega^1_{G^1_n}
\]
the homomorphism $\partial$ fits into the exact sequence:
\[
\xymatrix{ 0\ar[r] & \I/\I^2 \ar[r]^-\partial &
p_2^*\Omega^1_{G^1_n} \ar[r] & \Omega^1_{N^1_n/\cur_2}\ar[r]&0}
\]
Dualizing we obtain the exact sequence:
\[
\xymatrix{ 0\ar[r] &Hom(
\Omega^1_{N^1_n/\cur_2},\O_{N^1_n})\ar[r]&p_2^*T_{G^1_n}
\ar[r]^-{\partial^\vee}& N_{N^1_n} \ar[r] & T^1_{p_1} \ar[r] & 0}
\]
where $T^1_{p_1}$ is the first relative cotangent sheaf of $p_1$
and $N_{N^1_n}=N_{N^1_n(\cur)/\cur_2\times G^1_n}$. Since
$D_0(\partial)=D_0(\partial^{\vee})$,  we see that
$\Supp(D_0(\partial))= \Supp(T^1_{p_1})$, and this is precisely
the relative singular locus of $p_1$. By Lemma \ref{L:D0} for each
$x+y \in \cur_2$ we have
\[
\iota(\cur_2\times
 G^1_{n-2})\cap N^1_n(x+y) \subset \mathrm{Sing}(N^1_n(x+y))
\]
and therefore  $\iota(\cur_2\times
 G^1_{n-2})\subset \mathrm{Supp}(D_0(\partial))$.

(ii)  Note that $U$ is dense in $G^1_n$ (Lemma \ref{L:rho1}).
Since domain and codomain of $\partial$
  have ranks 1 and $\rho$ respectively,  it follows that
  every irreducible component of $F$ has codimension $\le \rho$ in $N^1_n(\cur)$,
  i.e. it has dimension  $\ge 1$.
Assume that $\cur$ is a general curve of genus $g$.   then a
general $(L,V)\in
 U$ has simple ramification (by a straightforward count of parameters), and therefore $s_{L,V}$ is
 nonsingular and irreducible, by Lemma \ref{L:sLV2})(ii). If
    we restrict $\partial^\vee$ to   $s_{L,V}= p_2^{-1}(L,V)$,  we  obtain:
\[
\xymatrix{ T_{(L,V)}G^1_n\otimes \O_{s_{L,V}} \ar[r] &
N_{s_{L,V}/\cur_2}}
\]
Since  $N_{s_{L,V}/\cur_2}$ has positive  degree $(n-1)^2-\gamma$
it follows that
 the vector bundle $p_2^*T_{G^1_n}^\vee\otimes N_{N^1_n}$ is $p_2$-relatively
 ample.
  We   now apply Theorem 0.3 of \cite{fS98} to deduce that $F \ne \emptyset$.
  This proves the non-emptiness statement in case $\cur$ is general. But the same property
  is preserved under specialization,
  and therefore it is also true only assuming that $\cur$ is Petri w.r. to pencils. \qed

\begin{prop}\label{P:NP2}
In the same hypothesis as before let $R_n=p_1(F) \subset \cur_2$.
Then:
\begin{description}
\item[(i)] $F$ is purely one-dimensional,   $R_n\subset \cur_2$ is
a curve and the projection $F \longrightarrow R_n$ is finite.
 \item[(ii)] For all $(x+y,(L,V))\in F$,   $x$ is not a
base point of the pencil $(L,V)$.
 \item[(iii)] If $\cur$ is general then  the curve $R_n$ does not contain the diagonal
 $\delta \subset \cur_2$.
\end{description}
\end{prop}

\proof (i)  It will suffice to prove the following:

\noindent (A) The restriction to $D_0(\partial)\cap p_2^{-1}(U)$ of the projection $p_1: N^1_n(\cur) \to \cur_2$ has no
positive dimensional fibres.

\noindent (B)  $p_1(D_0(\partial)\cap p_2^{-1}(U))$ does not contain the general point of the curve $x+\cur$ for
all $x \in \cur$.

Let's prove (A).
 Assume that for some $x+y \in \cur_2$ there is an irreducible curve
 $T\subset p_1^{-1}(x+y)\cap [D_0(\partial)\cap p_2^{-1}(U)] \subset N^1_n(x+y)\cap F$.
 Then,  for each $t\in T$ the corresponding   $(L_t,V_t)$ is base-point free and complete and, by Proposition \ref{P:NP1},
 is a singular point of $N^1_n(x+y)$. By the criterion (ii) of Corollary \ref{C:petri3},
  there is a unique effective divisor
  $D_t = z_1(t) + \cdots + z_{2\gamma-2n}(t)$ in $|\omega_\cur L_t^{-2}(x+y)|$.   Since  the divisors $D_t$ are not constant w.r. to $t \in T$, for some $t_0 \in T$ we have
  $x \in \Supp(D_{t_0})$. This implies that  $h^1(L_{t_0}^2(-y)) > 0$,  contradicting Proposition  \ref{P:kermu}(ii).
   This concludes the proof of (A).

  Let's prove (B). Assume that, given $x$,  for a general $y \in \cur$ there is
  \[
  (L_y,V_y) \in p_1^{-1}(x+y)\cap [D_0(\partial)\cap p_2^{-1}(U)]
  \]
   Then, by the criterion (ii) of Corollary \ref{C:petri3},
  for each such $y$ there is a unique effective divisor
  $D_y \in |\omega_\cur L_y^{-2}(x+y)|$.   Arguing as before we deduce
  that  $x \in \Supp(D_y)$ for some $y$ and we contradict  Proposition  \ref{P:kermu}(ii) again.

  (ii)  Clearly, if
$(x+y,(L,V))\in D_0(\partial)\cap p_2^{-1}(U)$ there is nothing to
prove;  moreover $h^1(L^2(-x-y))=1$, by Corollary \ref{C:petri3}.
Consider now $(x+y,(L,V))\in F$ arbitrary, and let $B$ be the
fixed divisor of $(L,V)$. By upper-semicontinuity we
  have $h^1(L^2(-x-y))>0$, and therefore $h^1(L^2(-x-y-B))>0$
as well. If $x < B$ then $h^1(L^2(-x-y-B))=0$, By Lemma
\ref{L:kerM}(c), and this is a contradiction.

(iii)    For a general $x \in \cur$ the curve $X = \cur \cup E$,
where $E$ is a general elliptic curve and $\cur \cap E= \{x\}$, is
a general element of $\Delta_1 \subset \overline{\M}_g$, where
$g=\gamma+1$. If $\delta \subset R_n$ then a pencil $(L,V) \in F
\cap p_1^{-1}(2x)$ descends to a pencil on $X$ whose Petri map is
not injective.  But from \cite{EH83} we know that the general
curve in $\Delta_1$ is Petri, and we get a contradiction.
   \qed

\section{Nodal curves}\label{S:nodcur}

Let $X$ be an integral projective curve of arithmetic genus $g \ge
4$ whose only singularity is an ordinary node $\bar z$. Let
\[
\nu: \xymatrix{ \cur \ar[r] & X}
\]
be the normalization, and let $\gamma:=g-1$ be its genus. Denote
by $x+y = \nu^{-1}(\bar z)$.   We have
$\nu^*\omega_X=\omega_\cur(x+y)$.

A torsion free rank 1 sheaf $L$ on $X$ is either invertible or it
is not invertible precisely at the point $\bar z$. In this case $L
= \nu_*\widetilde{L}$ where $\widetilde{L}$ is an invertible sheaf
on $\cur$. In particular the ideal sheaf $\I_{\bar z}$ of $\bar
z$, which is torsion free but not locally free, is
\[
\I_{\bar z} = \nu_*\O_\cur(-x-y)
\]

\begin{lemma}\label{L:cusp1}
 \begin{description}

 \item[(i)]
There is a canonical isomorphism of $\O_X$-modules
$$
Hom_X(\nu_*\O_\cur,\O_X) \cong \I_{\bar z}
$$
associating to $\varphi: \nu_*\O_\Gamma \longrightarrow \O_X$ the
section of $\O_X$ corresponding to the composition
$$
\xymatrix{ \O_X \ar[r]& \nu_*\O_\cur \ar[r]^-\varphi& \O_X}
$$

\item[(ii)] For every ideal sheaf ${\cal N} \subset \O_\cur$ we
have $\I_{\bar z}(\nu_*{\cal N}) \subset \O_X$ ($\I_{\bar z}$ is
the {\rm conductor} of $\O_\cur$ in $\O_X$).

\item[(iii)] Let ${\cal N} \subset \O_\cur$ be an ideal sheaf.
Then there is a canonical isomorphism:
$$
\nu_*[Hom_\cur({\cal N},\O_\cur)] \cong Hom_X(\I_{\bar
z}(\nu_*{\cal N}), \O_X)
$$

\end{description}
\end{lemma}

  \proof
  (i) is elementary and (ii) is true by definition.

  (iii)  is local around $\bar z$. So let  $O= \O_{X,\bar z}$, $A$ the integral closure of $O$,
  $M=\Hom_O(A,O)\subset O$ the maximal ideal,
  and
  $N={\cal N}_z$.  We have
  $$
  \begin{array}{clcc}
  \Hom_A(N,A) &= \Hom_A(MN,M)&&\hbox{(because $M$ is invertible)}  \\
  & =\Hom_A(MN,\Hom_O(A,O))  \\
  &=\Hom_O(MN,O)
  \end{array}
  $$
  The last equality is obtained by associating to $f: MN \to \Hom_O(A,O)$ the $O$-homomorphism
  $\varphi: MN \to O$ defined by $\varphi(a) = f(a)(1)$.
  \qed

 The {\it degree}  of a torsion free rank 1 sheaf $L$ on $X$ is
 $$
 \deg(L) = \chi(L) - \chi(\O_{X})
 $$

 \begin{lemma}\label{L:cusp1.5}
 If $L$ is not locally free, then $L=\nu_*\widetilde{L}$ for some invertible sheaf  $\widetilde{L}$ on $\cur$ and
 \begin{equation}\label{E:cusp1}
 \deg(\widetilde{L})=\deg(L)-1
 \end{equation}
 In particular
 $$
 \deg(\I_{\bar z}) = -1
 $$
 and
 $$
 \deg(\nu_*\O_\Gamma) = 1
 $$
 \end{lemma}

 \proof
 In fact  $\deg(\widetilde{L})= \chi(\widetilde{L}) - \chi(\O_\cur)$;
 we have $\chi(\widetilde{L}) = \chi(L)$ because $\nu$ is finite, and
 $$
 \chi(\O_\Gamma) = 1 - \gamma = 1 - g+1 = \chi(\O_{X})+1
 $$
 \qed

We need the following two definitions.

 \begin{definition}\label{D:cusp1}
 A $g^r_n$  on $X$ is a pair $(L,V)$ where $L$ is an invertible sheaf of degree $n$ and
$V \subset H^0(L)$ is a subspace of dimension  $r+1$.
 A \emph{ generalized $g^r_n$} on $X$ is a pair
$(L,V)$, where $L$ is a torsion free rank 1 sheaf of degree $n$
and $V \subset H^0(L)$ is a vector subspace of dimension $r+1$.
Obviously, every $g^r_n$ is also a generalized $g^r_n$.
\end{definition}

 \begin{definition}\label{D:cusp2}
  Let $(L,V)$ be a generalized $g^r_n$ on $X$. The natural map
$$
\mu_0(V): \xymatrix{V \otimes \Hom(L,\omega_{X}) \ar[r]&
H^0(\omega_{X})}
$$
is called the \emph{Petri map of $(L,V)$}. If $L$ is invertible
then $\mu_0(V)$ is just the ordinary Petri map of $(L,V)$:
 \[
 \mu_0(V): \xymatrix{
  V \otimes H^0(\omega_XL^{-1})\ar[r] & H^0(\omega_X) }
  \]
   If
$\mu_0(V)$ is injective then $(L,V)$ is called {\rm a Petri
(generalized) $g^r_n$}.
 Given $r,n$ the curve $X$   is
called  a {\rm Petri curve w.r. to generalized $g^r_n$'s}  if
every generalized $g^r_n$ on $X$ is Petri. Similarly, given $r \ge
0$, we call $X$ a {\rm Petri curve w.r. to generalized $g^r$'s} if
for each $n$ every limit $g^r_n$ on $X$ is Petri.
 We call $X$ a {\rm Petri curve} if it is Petri w.r. to generalized $g^r_n$'s for all $r,n$.
 \end{definition}

If $L$ is a  torsion free rank 1 sheaf of degree $n$ on $X$ which
is not locally free, then $L= \nu_*\widetilde{L}$ where
 $\widetilde{L}$ is a line bundle of degree $n-1$ on $\cur$;  we have
$$
H^0(X,L) = H^0(\cur,\widetilde{L})
$$
 and
 $$
 \begin{array}{clr}
Hom_X(L,\omega_{X}) &= Hom_X(\nu_*\widetilde{L},\omega_{X})\\
&=Hom_X(\I_{\bar z}[\nu_*\widetilde{L}(x+y)],\omega_{X})\\
&= \nu_*[Hom_\cur(\widetilde{L}(x+y),\nu^*\omega_X)] & \hbox{(by Lemma \ref{L:cusp1}(iii))} \\
&= \nu_*[Hom_\cur(\widetilde{L},\omega_\cur)]
\end{array}
$$
so that
\begin{equation}\label{E:cusp4}
\Hom_X(L,\omega_{X}) = \Hom_\cur(\widetilde{L},\omega_\cur)=
H^0(\cur,\omega_\cur \widetilde{L}^{-1})
 \end{equation}
 If $L$ is invertible on $X$ then we have an inclusion:
 \begin{equation}\label{E:cusp2}
 H^0(X,L) \subset H^0(X,\nu_*\nu^*L) = H^0(\cur,\nu^*L)
 \end{equation}
After these preliminaries we can now prove a proposition which
relates the Petri maps on $X$ with certain maps on $\cur$.

\begin{prop}\label{P:cusp2}
Let $(L,V)$ be a generalized $g^r_n$ on $X$. Then
\begin{description}
\item[(i)] If $L$ is invertible we have the commutative diagram:
$$
\xymatrix{ V \otimes H^0(\omega_XL^{-1})\ar@{^(->}[d] \ar[r]^-{\mu_0(V)}&  H^0(\omega_{X})\ar@{=}[d]  \\
  V\otimes H^0(\omega_\cur \nu^*L^{-1}(x+y))\ar[r]^-{M}
   & H^0(\omega_\cur(x+y))}
 $$
where $M$ is the natural multiplication map. If moreover $\bar z$
is a base point of $(L,V)$ then we have a commutative diagram:
$$
\xymatrix{ V \otimes H^0(\omega_XL^{-1})\ar@{^(->}[d] \ar[r]^-{\mu_0(V)}&  H^0(\omega_{X})   \\
  V\otimes H^0(\omega_\cur \nu^*L^{-1}(x+y))\ar[r]^-{M}
   & H^0(\omega_\cur)\ar@{^(->}[u]}
 $$

\item[(ii)]  If $L$ is not invertible then $L= \nu_*\widetilde{L}$
for some invertible sheaf $\widetilde{L}$ on $\Gamma$ and we have
a commutative diagram:
$$
\xymatrix{
V \otimes \Hom(L,\omega_{X}) \ar[r]^-{\mu_0(V)}\ar@{=}[d]&  H^0(\omega_{X}) \\
 V\otimes H^0(\omega_\cur \widetilde{L}^{-1}) \ar[r] &
 H^0(\omega_\cur)\ar@{^(->}[u]}
 $$

 \end{description}

\end{prop}

\proof (i) In the first diagram the left vertical inclusion is
induced by the inclusion (\ref{E:cusp2}) applied to
$\omega_XL^{-1}$. The commutativity is clear. The second diagram
holds because $V= V(-x-y)$ since $\bar z$ pulls back to $x+y$.

(ii)  The left  vertical equality is a consequence of
(\ref{E:cusp4}) and the commutativity is clear in this case too.
\qed

\begin{corollary}\label{C:nodal1}
   If      $\cur$ is Petri w.r. to pencils
and  the map
\begin{equation}\label{E:cusp3}
 \xymatrix{M: \nu^*V\otimes H^0(\omega_\cur
(\nu^*L)^{-1}(x+y))\ar[r]& H^0(\omega_\cur(x+y))}
\end{equation}
is injective for all $g^1_n$'s    on
 $\Gamma$  of the form $(\nu^*L,\nu^*V)$ for which $x$ is not
 a base point then $X$ is Petri w.r. to generalized $g^1_n$'s.
\end{corollary}

\proof   If $\cur$ is Petri  w.r. to pencils, then the map
(\ref{E:cusp3}) is injective whenever $x$ is a base point of
$(\nu^*L,\nu^*V)$, because $(\nu^*L,\nu^*V)\in N^1_n(x+y)$ and by
Lemma \ref{L:kerM}(c). Therefore the corollary is an immediate
consequence of Proposition \ref{P:cusp2}, since, by
\ref{P:cusp2}(ii),   the
 Petri map is  injective for generalized $g^1_n$'s $(L,V)$ for
 which $L$ is not locally free.  \qed

 As a consequence of all the above we have:

 \begin{theorem}\label{T:nodal1}
Assume that  $\cur$ is  Petri w.r. to pencils and let $R_n \subset
\cur_2$  be the curve introduced in  Proposition \ref{P:NP2}.
Then:
\begin{description}
\item[(i)] If $x+y\in \cur_2\setminus R_n$ then $X$ is Petri w.r.
to generalized $g^1_n$'s.
 \item[(ii)] If $x+y \in R_n$ then $X$
has finitely many  generalized $g^1_n$'s $(L,V)$ for which the
Petri map $\mu_0(V)$ is not injective, and for all of them  $L$ is
invertible (i.e. they are $g^1_n$'s) and $x$ is not a base point
of $(L,V)$. If moreover $x+y$ is general in $R_n$ then all the
non-Petri $g^1_n$'s on $X$ are complete.

\item[(iii)]  If $x+y \in \cur_2$ is general then $X$ is Petri
w.r. to generalized pencils.
\end{description}
\end{theorem}

\proof (i) If we choose $x+y\notin R_n$
 then for all  $(\nu^*L,\nu^*V)$ for which $x$ is not a base point
 and $L$ is invertible the   map  (\ref{E:cusp3}) is injective, by
 Propositions \ref{P:TN} and \ref{P:NP1}. Therefore we can apply
 Corollary \ref{C:nodal1}.

 (ii) by
\ref{P:cusp2}(ii) and Corollary \ref{C:nodal1},   the
 Petri map is  injective for generalized $g^1_n$'s $(L,V)$ for
 which $L$ is not locally free or $\bar z$ is a base point.
 Let $(L,V)$ be a generalized $g^1_n$ on $X$ for which $\bar z$ is not a base point
 and $L$ is invertible. If the map (\ref{E:cusp3}) is
 not injective then it follows from
 Propositions \ref{P:TN} and \ref{P:NP1} that
 $(\nu^*L,\nu^*V)\in F$. Since the map
 $F \longrightarrow p_1(F)$ is finite, we deduce that there are finitely many
 $g^1_n$'s on $X$   such that the map (\ref{E:cusp3}) is
 not injective.  The assertion about the completeness is a consequence of the definition of $F$.

 (iii) If $x+y$ is general then it is not on the finitely many
 curves $R_n$,  $\frac{1}{2}g \le n \le g$, and therefore $X$ is
 Petri w.r. to generalized $g^1_n$'s for all such $n$'s, by (i)
 above, and therefore it is Petri w.r. to generalized pencils.
 \qed

\section{Proof of Theorem \ref{T:main2}}\label{S:main2}

\proof   Let $X$ be as in Theorem \ref{T:nodal1}(ii), with $x+y
\in R_n$ general. Then $X$ is irreducible and 1-nodal of
arithmetic genus $g$, has a non-zero finite number of generalized
$g^1_n$'s $(L,V)$ which are non-Petri, for all of them $L$ is
invertible, they are base-point free and complete.  We consider a
modular family of deformations of $X$:
\[\xymatrix{
X\ar@{^(->}[r]\ar[d] & \X \ar[d]^-f \\
\mathrm{Spec}(\mathbb{C})\ar[r]^-{b_0} & B}
\]
In other words $B$ is nonsingular of dimension $3g-3$ and  the
Kodaira-Spencer map
\[\xymatrix{
\kappa_b: T_bB
\ar[r]&\mathrm{Ext}^1(\Omega^1_{\X(b)},\O_{\X(b)})}\] is an
isomorphism for all $b \in B$.  We will further assume that all
fibres of $f$ are in $\M_g \cup \Delta_0$ and that $f$ has a
section $\sigma: B \longrightarrow \X$. We have a commutative
diagram:
\begin{equation}\label{E:Jn1}
 \xymatrix{
J_n(\X/B)\times_B \X \ar[r]^-{p_2} \ar[d]_-{p_1}&
 \X \ar[d]^-f \\
J_n(\X/B)\ar[r]_-q &B}
 \end{equation}
where $J_n(\X/B)$ is the relative Picard variety parametrizing
invertible sheaves of degree $n$ on the fibres of $f$. We now
consider the relative scheme $G^1_n(\X/B)$ of $g^1_n$'s on the
fibres of $f$. It can be defined as in the (absolute) case of a
single curve, using the existence of a relative Poincar\'e sheaf
$\cP$ on $J_n(\X/B)\times_B \X$. We need to recall  its
construction. One fixes an integer $m \ge 2g-1-n$ and considers
the relative grassmannian
\[\pi:G(2,p_{1*}[\cP\otimes p_2^*\O_\X(m\sigma(B))])\longrightarrow J_n(\X/B)\] and the
tautological locally free sheaf of rank two $\mathcal{E}
\subset\pi^*p_{1*}[\cP\otimes p_2^*\O(m\sigma(B))]$. Then the
evaluation homomorphism:
\[\xymatrix{
\cP\otimes p_2^*\O_\X(m\sigma(B))\ar[r] & \cP\otimes
p_2^*\O_{m\sigma(B)}(m\sigma(B))}\] induces a homomorphism:
\[\xymatrix{
\pi^*p_{1*}[\cP\otimes p_2^*\O_\X(m\sigma(B))]\ar[r]&
\pi^*p_{1*}[\cP\otimes p_2^*\O_{m\sigma(B)}(m\sigma(B))]}\] and
$G^1_n(\X/B)$ is the vanishing scheme of the composition:
\[
\xymatrix{\mathcal{E} \ar[r]&\pi^*p_{1*}[\cP\otimes
p_2^*\O_\X(m\sigma(B)) ]\ar[r]& \pi^*p_{1*}[\cP\otimes
p_2^*\O_{m\sigma(B)}(m\sigma(B))]}\] The restriction
$\mathcal{E}_{|G^1_n(\X/B)}$  will be denoted   by
$\widetilde{\mathcal{E}}$.

 In a similar way one defines the relative scheme
$G^{g-n}_{2g-2-n}(\X/B)$ as a closed subscheme of the relative
grassmannian
\[
\chi:G(g-n+1,p_{1*}[\cP^{-1}\otimes
p_2^*\omega_{\X/B}(e\sigma(B))])\longrightarrow J_{2g-2-n}(\X/B)
\]
for some $e \ge n+1$, and the restricted tautological sheaf
$\widetilde{\mathcal{F}}$ of rank $g-n+1$ on
$G^{g-n}_{2g-2-n}(\X/B)$. We denote by
\[
 \overline{\chi}: G^{g-n}_{2g-2-n}(\X/B) \longrightarrow J_{n}(\X/B)
\]
the composition of $\chi$ with the isomorphism
$J_{2g-2-n}(\X/B)\longrightarrow J_{n}(\X/B)$ induced by
residuation with respect to $\omega_{\X/B}$.
 Now consider the
following diagram:
\[
\xymatrix{
\widetilde{G}\ar[r]^-{q_2}\ar[d]_-{q_1}&G^{g-n}_{2g-2-n}(\X/B)\ar[d]^-{\overline{\chi}
}\\
G^1_n(\X/B) \ar[r]_\pi & J_{n}(\X/B)\ar[r]_-q&B}
\]
where we denoted by $\widetilde{G} =
G^1_n(\X/B)\times_{J_{n}(\X/B)}G^{g-n}_{2g-2-n}(\X/B)$.  Let
$\omega := (q\pi q_1)^*f_*(\omega_{\X/B})$. It is locally free of
rank $g$ on $\widetilde{G}$, and we have a natural multiplication
homomorphism:
\[\xymatrix{
\widetilde{\mu_0}:q_1^*\widetilde{\mathcal{E}}\otimes
q_2^*\widetilde{\mathcal{F}}\ar[r] & \omega}\]
Consider the
degeneracy locus:
\[
P^1_{g,n}(\X/B) := D_{2(g-n+1)-1}(\widetilde{\mu_0}) \subset
\widetilde{G}\] It is non-empty because its projection to
$G^1_n(\X/B)$ contains every $g^1_n$ $(L,V)$ on $X$ such that
$\mu_0(L,V)$ is not injective. Moreover $P^1_{g,n}(\X/B)$ has
expected codimension
\[g-2(g-n+r)+1 = \rho(g,1,n)+1\]
in $\widetilde{G}$.  Therefore every irreducible component
$\widetilde{P}$ of $P^1_{g,n}(\X/B)$ satisfies:
\begin{equation}\label{E:dimP}
    \dim(\widetilde{P}) \ge \dim(\widetilde{G})-[\rho(g,1,n)+1]\ge 3g-4
\end{equation}
Consider in particular a component $\widetilde{P}$ such that
$q_1(\widetilde{P})$ contains an $(L,V)$ on $X$ with $\mu_0(L,V)$
not injective. Let
\[
P:=(q\pi q_1)(\widetilde{P}) \subset B
\]
Then $b_0\in P$ and the induced morphism
$\widetilde{P}\longrightarrow P$ is finite above $b_0$,   since
$q_1$ is birational above $(L,V)$ (because $(L,V)$ is complete)
and $q_1(\widetilde{P})\cap G^1_n(\X/B)(b_0)$ is finite. Therefore
$\dim(P)= \dim(\widetilde{P}) \ge 3g-4$. But $\dim(P) \le 3g-4$ as
well because $P$ is a closed locus parametrizing  curves that are
non-Petri w.r. to pencils, and therefore it does not contain all
the points of $B$ parametrizing curves in $\Delta_0$, by Theorem
\ref{T:nodal1}(iii).
 Let $\beta: B \longrightarrow  \M_g\cup \Delta_0$ be the
 functorial morphism, which is finite by construction. Then
 $\beta(P)$ is an irreducible divisor containing $[X]$, but not entirely contained in $\Delta_0$
 (Theorem \ref{T:nodal1}) and
 $\beta(P)\cap \M_g$ is a divisorial component of $P^1_{g,n}$
 having the required properties.  This concludes the proof of
 Theorem \ref{T:main2}.

\noindent\textsc{address of the authors:}

\noindent Dipartimento di Matematica,
  Universit\`a Roma Tre \newline Largo S. L. Murialdo 1,
  00146 Roma, Italy.
  \smallskip
\newline   \texttt{bruno@mat.uniroma3.it}
\newline  \texttt{sernesi@mat.uniroma3.it}

\end{document}